\documentclass{ifacconf}

\usepackage{graphicx}      % include this line if your document contains figures
\usepackage{natbib}        % required for bibliography

\usepackage{amsmath,amssymb,amsfonts}
\usepackage{enumerate}
\usepackage[svgnames]{xcolor}
\usepackage{tikz}
\graphicspath{ {img/} }

\newtheorem{casest}{Case Study}

\usetikzlibrary{arrows}

\tikzset{
	treenode/.style = {align=center, inner sep=0pt, text centered,
		font=\sffamily},
	task/.style = {treenode, circle, white, font=\sffamily\bfseries, draw=black,
		fill=black, text width=0.9em},% arbre rouge noir, noeud noir
	agent/.style = {treenode, circle, SteelBlue, draw=SteelBlue,
		text width=0.9em, thin},% arbre rouge noir, noeud rouge
}

\begin{document}
\begin{frontmatter}

\title{Exploiting Structure in the Bottleneck Assignment Problem\thanksref{footnoteinfo}} 

\thanks[footnoteinfo]{The research is funded by Defence Science and Technology Group through research agreements MyIP: 7558 and MyIP: 7562.}

\author[First]{Mitchell Khoo} 
\author[Second]{Tony A. Wood} 
\author[Second]{Chris Manzie}
\author[Second]{Iman Shames}

\address[First]{Department of Electrical and Electronic Engineering at the University of Melbourne, (e-mail: khoom1@student.unimeld.edu.au).}
\address[Second]{(e-mail: \{wood.t,manziec,iman.shames\}@unimelb.edu.au)}
%\address[Third]{(e-mail: manziec@unimelb.edu.au)}
%\address[Fourth]{(e-mail: iman.shames@unimelb.edu.au)}

\begin{abstract}                % Abstract of not more than 250 words.
An assignment problem arises when there exists a set of tasks that must be allocated to a set of agents. The bottleneck assignment problem (BAP) has the objective of minimising the most costly allocation of a task to an agent. Under certain conditions the structure of the BAP can be exploited such that subgroups of tasks are assigned separately with lower complexity and then merged to form a combined assignment. In particular, we discuss merging the assignments from two separate BAPs and use the solution of the subproblems to bound the solution of the combined problem. We also provide conditions for cases where the solution of the subproblems produces an exact solution to the BAP over the combined problem. We then introduce a particular algorithm for solving the BAP that takes advantage of this insight. The methods are demonstrated in a numerical case study.
\end{abstract}

\begin{keyword}
Algorithms; Decision-making; Graph theory; Agents; Optimization problems.
\end{keyword}

\end{frontmatter}
%===============================================================================

\section{Introduction}

An assignment problem arises when multiple tasks are to be allocated to multiple agents. For example, situations where jobs are to be assigned to a group of workers or passengers positioned at different locations in a city are to be picked up by a fleet of cars. Tasks can be assigned based on many different criteria. See \cite{taxonomy}, \cite{assignprob}, and \cite{survey} for reviews on the different objectives for assignment problems.

One particular objective is to assign tasks to agents such that the total cost of the assignment is minimised. This type of assignment problem is called the linear assignment problem (LAP). The Hungarian Method in \cite{hungarian} is a well-studied algorithm for solving the LAP. In \cite{distHungarian}, a distributed version of the Hungarian Method is presented; a distributed algorithm is one that does not rely on a centralised decision-maker for computation. In \cite{auc1} and \cite{auc2}, so-called auction algorithms are presented to solve the LAP. A greedy algorithm is one where tasks are allocated to agents sequentially. Each allocation is made according to the lowest cost amongst the remaining choices. In \cite{CBAA}, the Consensus-Based Auction Algorithm (CBAA) is presented, which is a greedy algorithm used to obtain suboptimal solutions to the LAP with low computational cost compared to algorithms for solving the LAP exactly.

Another objective is to assign tasks to agents such that the costliest allocation is minimised, which corresponds to the bottleneck assignment problem (BAP). The BAP has application in time-critical problems. For example in \cite{decoys}, a set of decoys must travel to a set of positions such that the worst-case positioning time is minimised. In \cite{BAP1}, a threshold algorithm is presented, where a threshold is iteratively increased until it is possible to find an assignment containing only allocations of tasks to agents with costs smaller than the threshold. In \cite{BAP3,BAP4}, the bound on the completion time of the threshold algorithm is reduced moving the threshold according to a binary search pattern. In \cite{BAP2}, an algorithm is presented that iteratively solves the BAP over an increasing subset of agents and tasks. The subset size is increased until it contains all the agents and tasks. In \cite{distBAP}, a distributed algorithm for solving the BAP is introduced. There are other variants of the BAP. Such variants include the scheduling problems in \cite{bus} and \cite{categorisation}, which require assigning more than one task per agent. In fact, this can be regarded as an example of a time-extended assignment from the taxonomy in \cite{taxonomy}.

In this paper, we focus on the BAP and restrict the scope to having each agent carry out at most one task and each task requiring at most one agent for completion. The contribution of this work is to investigate structure that can be exploited to solve the BAP efficiently. Consider partitioning the sets of agents and tasks, i.e., splitting the assignment problem into two smaller BAPs. We can use the two solutions of the subproblems for solving the combined BAP. Consider the following three ways to exploit the structure of the BAP. We relate each scenario to a ride-sharing application for illustration.

For the first scenario, assume the sets of agents and tasks were partitioned equitably, i.e., none of the subproblems has fewer tasks than agents. Merging the solutions of the two subproblems forms a valid but possibly suboptimal assignment in the combined problem. In fact, the cost of the merged assignment is an upper bound on the cost of the optimal BAP solution. In a ride-sharing application, two rival companies may assign their own vehicles to their own customers. However, they may find that pooling their resources allows a better service for all customers.

For the second scenario, we define a bottleneck cluster as a group of agents and tasks with small allocation costs amongst each other. When the two subproblems consist of two separate bottleneck clusters, we can determine conditions under which the solutions of the subproblems form an exact solution to the combined problem. Consider two cities each with their own sets of vehicles and customers. If the cities were geographically far apart, there is no benefit for vehicles in one city to serve customers in the other city.

The third scenario relates to the algorithm in \cite{distBAP}. Knowing the solutions to the two subproblems leads to information about task-to-agents allocations that are particularly costly. We can eliminate suboptimal options when the algorithm is initialised to solve the combined problem. Assume a group of customers has been assigned to vehicles. Then, a new group of customers requests to be picked up. By only considering the idle vehicles for the new customers and not the previously assigned ones, the resulting assignment problem has lower complexity. The solution from the two subproblems can be used as a warm-start to solving the combined problem.

\section{Preliminaries} \label{prelim}

Given an arbitrary undirected graph $\mathcal{G}=(\mathcal{V,E})$ with vertex set $\mathcal{V}$ and edge set $\mathcal{E}$, consider the following definitions found in \cite{hk_alg} and \cite{distBAP}.

\begin{defn}[Maximum Cardinality Matching] \label{def:matching}
	A matching $\mathcal{M}$ in $\mathcal{G}$ is a set of edges such that $\mathcal{M}\subseteq \mathcal{E}$ and no vertex $v\in\mathcal{V}$ is incident with more than one edge in $\mathcal{M}$. A maximum cardinality matching (MCM) is a matching $\mathcal{M}_{max}$ in $\mathcal{G}$ of maximum cardinality.
\end{defn}

Let $\mathcal{A}_b$ be a set of agents and $\mathcal{B}_b$ be a set of tasks, where $\mathcal{A}_b\cap \mathcal{B}_b=\emptyset$. Consider an arbitrary complete bipartite graph $\mathcal{G}_b=(\mathcal{V}_b,\mathcal{E}_b)$ with vertex set $\mathcal{V}_b=\mathcal{A}_b\cup \mathcal{B}_b$ and edge set $\mathcal{E}_b=\{\{i,j\}|i\in\mathcal{A}_b,j\in\mathcal{B}_b\}$. Let $\mathcal{C}(\mathcal{G}_b)$ be the set of all MCMs of $\mathcal{G}_b$. Let $w:\mathcal{E}_b\mapsto \mathbb{R}$ map edges to real-valued weights. The BAP for graph $\mathcal{G}_b$ is formulated as
\begin{align} \label{eq:BAP}
BOT(\mathcal{G}_b):\qquad&\min_{\mathcal{M}\in \mathcal{C}(\mathcal{G}_b)}\max_{\{i,j\}\in \mathcal{M}} \quad w(\{i,j\}).
\end{align}

\begin{defn}[Bottleneck edge] \label{def:botedge}
	A bottleneck edge of graph $\mathcal{G}_b$ is any $e\in \arg\max_{\{i,j\}\in\mathcal{M}}w(\{i,j\})$, for any MCM $\mathcal{M}$ that is a minimiser of $BOT(\mathcal{G}_b)$.
\end{defn}

\begin{defn}[Neighbours] \label{def:neigh}
	The set of neighbours of vertex $v\in\mathcal{V}$ in $\mathcal{G}$ is defined as $N_v=\{k|\{v,k\}\in \mathcal{E}\}$.
\end{defn}

Note that given a vertex $v\in\mathcal{V}$, $\forall k\in N_v$, $v\in N_k$.

\begin{defn}[Path] \label{def:path}
	Let a sequence of distinct vertices $v_1,v_2,...,v_{l+1}\in \mathcal{V}$ be such that for $k=1,2,...,l$, $v_{k+1}\in N_{v_{k}}$. The set of edges $\mathcal{P}=\{ \{v_k,v_{k+1}\}\}_{k=1,2,...,l}$ is then said to be a path between $v_1$ and $v_{l+1}$, with length $l$.
\end{defn}

\begin{defn}[Alternating path] \label{def:alt}
	Given a matching $\mathcal{M}$ and a path $\mathcal{P}$, $\mathcal{P}$ is an alternating path relative to $\mathcal{M}$ if and only if each vertex that is incident to an edge in $\mathcal{P}$ is incident with no more than one edge in $\mathcal{P}\cap\mathcal{M}$ and no more than one edge in $\mathcal{P}\backslash\mathcal{M}$.
\end{defn}

\begin{defn}[Free vertex] \label{def:free}
	Given a matching $\mathcal{M}$, a vertex $v\in \mathcal{V}$ is free if and only if for all $w\in \mathcal{V}$, $\{v,w\}\notin \mathcal{M}$.
\end{defn}

\begin{defn}[Augmenting path] \label{def:aug}
	Given a matching $\mathcal{M}$ and a path $\mathcal{P}$ between vertices $v$ and $v'$, $\mathcal{P}$ is an augmenting path relative to $\mathcal{M}$ if and only if $\mathcal{P}$ is an alternating path relative to $\mathcal{M}$ and $v$ and $v'$ are both free vertices.
\end{defn}

\begin{defn}
	(Alternating tree) Given a matching $\mathcal{M}$, $\mathcal{G}$ is an alternating tree relative to $\mathcal{M}$ if and only if $\mathcal{G}$ is a tree and any path between the root vertex of $\mathcal{G}$ and every other vertex in $\mathcal{G}$ is an alternating path relative to $\mathcal{M}$.
\end{defn}

\section{Problem Formulation} \label{formulation}

Let there be two sets of agents $\mathcal{A}_1=\{a_1,a_2,...,a_{m_1}\}$ and $\mathcal{A}_2=\{\alpha_1,\alpha_2,...,\alpha_{m_2}\}$ and two sets of tasks $\mathcal{B}_1=\{b_1,b_2,...,b_{n_1}\}$ and $\mathcal{B}_2=\{\beta_1,\beta_2,...,\beta_{n_2}\}$. Define the sets $\mathcal{A}_3:=\mathcal{A}_1\cup\mathcal{A}_2$ and $\mathcal{B}_3:=\mathcal{B}_1\cup\mathcal{B}_2$. Let $m_3=m_1+m_2$ and $n_3=n_1+n_2$ and assume $m_1\geq n_1$ and $m_2\geq n_2$. For $k=1,2,3$, define $\mathcal{V}_k:=\mathcal{A}_k\cup \mathcal{B}_k$, $\mathcal{E}_k:=\{\{i,j\}|i\in\mathcal{A}_k,j\in\mathcal{B}_k\}$ and graph $\mathcal{G}_k:=(\mathcal{V}_k,\mathcal{E}_k)$. Define $\mathcal{D}(\mathcal{G}_b):=\arg\min_{\mathcal{M}\in \mathcal{C}(\mathcal{G}_b)}\max_{\{i,j\}\in \mathcal{M}}w(\{i,j\})$, the set of solutions to $BOT(\mathcal{G}_b)$ for any bipartite graph $\mathcal{G}_b$.

\begin{assum} \label{as:bot1}
	Assume we have $\mathcal{M}_1\in \mathcal{D}(\mathcal{G}_1)$ and $e_1\in\arg\max_{\{i,j\}\in\mathcal{M}_1}w(\{i,j\})$, i.e., an arbitrary solution to $BOT(\mathcal{G}_1)$ and a corresponding bottleneck edge of $\mathcal{G}_1$.
\end{assum}

\begin{assum} \label{as:bot2}
	Assume we have $\mathcal{M}_2\in \mathcal{D}(\mathcal{G}_2)$ and $e_2\in\arg\max_{\{i,j\}\in\mathcal{M}_2}w(\{i,j\})$, i.e., an arbitrary solution to $BOT(\mathcal{G}_2)$ and a corresponding bottleneck edge of $\mathcal{G}_2$.
\end{assum}

\begin{prob} \label{prob:sol1}
	Given Assumptions \ref{as:bot1} and \ref{as:bot2}, find a solution to $BOT(\mathcal{G}_3)$, i.e., find some matching $\mathcal{M}_3\in \mathcal{D}(\mathcal{G}_3)$.
\end{prob}

%Given a matching $\mathcal{M}_3$ defined in Problem~\ref{prob:sol1}, a bottleneck edge of $\mathcal{G}_3$ is $e_3\in\arg\max\limits_{\{i,j\}\in\mathcal{M}_3}w(\{i,j\})$.
In Section~\ref{sec:structure}, we define structures of the BAP that can be exploited to solve Problem~\ref{prob:sol1}. Then in Section~\ref{section:distBAP}, we discuss a specific algorithm that allows us to exploit some structure of the BAP discussed in Section~\ref{sec:structure}.

\section{Structure of the BAP} \label{sec:structure}

In this section, we discuss structures of the BAP that can be exploited. We first introduce an upper bound on the weight of a bottleneck edge of $\mathcal{G}_3$, in terms of the bottleneck edges $\mathcal{G}_1$ and $\mathcal{G}_2$. Then, we introduce bottleneck clusters and provide conditions when the solution to Problem~\ref{prob:sol1} is found by merging matchings $\mathcal{M}_1$ and $\mathcal{M}_2$.

\subsection{A Bound on the Optimal BAP Solution} \label{subsec:bound}

\begin{thm} \label{thm:bound}
	Under Assumptions \ref{as:bot1} and \ref{as:bot2}, it holds that $\max\{w(e_1),w(e_2)\}$ is an upper bound on $w(e_3)$, where $e_3$ is a bottleneck edge of $\mathcal{G}_3$.
\end{thm}
\begin{pf}
	By definition, $w(e_3)=\max_{\{i,j\}\in\mathcal{M}_3}w(\{i,j\})\leq \max_{\{i,j\}\in\mathcal{M}}w(\{i,j\})$ for any arbitrary $\mathcal{M}\in\mathcal{C}_3$. Since $\tilde{\mathcal{M}}=\mathcal{M}_1\cup\mathcal{M}_2\in\mathcal{C}_3$, $w(e_3)\leq \max_{\{i,j\}\in\tilde{\mathcal{M}}}w(\{i,j\})=\max\{w(e_1),w(e_2)\}$.\qed
\end{pf}

Given Assumptions~\ref{as:bot1} and~\ref{as:bot2}, the set $\mathcal{M}_1\cup\mathcal{M}_2$ is an MCM of $\mathcal{G}_3$. This MCM is possibly suboptimal to $BOT(\mathcal{G}_3)$. However, this bound of the BAP allows us to make a decision before solving Problem~\ref{prob:sol1} exactly. If the suboptimal solution is sufficient in our application, there is no need to invest further resources to solve Problem~\ref{prob:sol1} exactly.

%\begin{prob} \label{prob:condition}
%	Given Assumptions \ref{as:bot1} and \ref{as:bot2}, find conditions under which $\mathcal{M}_3=\mathcal{M}_1\cup\mathcal{M}_2$ is an exact solution to Problem~\ref{prob:sol1}.
%\end{prob}

%For solving the BAP, an augmenting path in Definition~\ref{def:aug} is a commonly used tool. We introduce the novel concept of a bottleneck cluster in the following subsection. Combining the concepts of augmenting paths and bottleneck clusters will allow us to exploit the structure of the BAP and provide conditions under which $\mathcal{M}_3=\mathcal{M}_1\cup\mathcal{M}_2$ is an exact solution to Problem~\ref{prob:sol1}.

\subsection{Bottleneck Clusters} \label{subsec:cluster}

We now introduce the novel concept of a bottleneck cluster. In \cite{distBAP}, conditions for determining if an edge is a bottleneck edge of a given graph are presented. We build on this result and discuss corresponding conditions under which $\mathcal{M}_3=\mathcal{M}_1\cup\mathcal{M}_2$ is an exact solution to $BOT(\mathcal{G}_3)$ when $\mathcal{G}_1$ and $\mathcal{G}_2$ are both bottleneck clusters.

Once again, consider an arbitrary complete bipartite graph $\mathcal{G}_b=(\mathcal{V}_b,\mathcal{E}_b)$. Given $\mathcal{M}$ is an MCM of $\mathcal{G}_b$, we define $\phi(\mathcal{G}_b,\mathcal{M}):=\mathcal{M} \cup \{e\in\mathcal{E}_b|w(e)<\max_{e'\in\mathcal{M}} w(e')\}$, the union of $\mathcal{M}$ and the set of edges that have weight strictly smaller than the largest edge in $\mathcal{M}$. With this tool we define a bottleneck cluster and a critical bottleneck edge.

\begin{defn}[Bottleneck cluster] \label{def:cluster}
	Let $\mathcal{M}_b$ be a solution to $BOT(\mathcal{G}_b)$. Let $e=\{a_b,b_b\}\in\mathcal{M}_b$ be a bottleneck edge of $\mathcal{G}_b$. Graph $\mathcal{G}_b$ is a bottleneck cluster relative to $e$ if and only if for any vertex $v\in\mathcal{V}_b$, there exists an alternating path $P$ between $v$ and bottleneck task $b_b$ such that $P\subseteq\phi(\mathcal{G}_b,\mathcal{M}_b)$.
\end{defn}

\begin{defn}[Critical bottleneck edge] \label{def:crit}
	Let $\mathcal{M}$ be an MCM of graph $\mathcal{G}_b$. Edge $e_c$ is a critical bottleneck edge of $\mathcal{G}_b$ relative to $\mathcal{M}$ if and only if $e_c\in\arg\max_{e\in\mathcal{M}} w(e)$ and $\phi(\mathcal{G}_b,\mathcal{M})\backslash\{e_c\}$ does not contain an augmenting path relative to $\mathcal{M}\backslash \{e_c\}$.
\end{defn}

Lemma~\ref{lem:distBAP} allows us to find a new MCM, which will have at least one less edge with weight $\max_{e\in\mathcal{M}} w(e)$ than $\mathcal{M}$.

\begin{lem}[Proof in \cite{distBAP}] \label{lem:distBAP}
	Let $\mathcal{M}$ be an arbitrary MCM of graph $\mathcal{G}_b$. Consider an edge $e\in\arg\max_{e\in\mathcal{M}} w(e)$. An augmenting path $P\subseteq \phi(\mathcal{G}_b,\mathcal{M})\backslash\{e\}$ exists relative to $\mathcal{M}\backslash\{e\}$ if and only if there exists an MCM $\mathcal{M}'$ of $\mathcal{G}_b$ such that $\mathcal{M}'\subseteq \phi(\mathcal{G}_b,\mathcal{M})\backslash\{e\}$.
\end{lem}

\begin{cor} \label{cor:crit}
	From Lemma~\ref{lem:distBAP}, it follows that every critical bottleneck edge of $\mathcal{G}_b$ is a bottleneck edge of $\mathcal{G}_b$.
\end{cor}

An MCM that is a solution to $BOT(\mathcal{G}_b)$ may contain more than one critical bottleneck edge. The following proposition shows how a critical bottleneck edge forms a particular alternating path between the bottleneck agent and bottleneck task.

\begin{assum} \label{as:botedge}
	Let $\mathcal{M}_b$ be a solution to $BOT(\mathcal{G}_b)$ and let $e_c=\{a_c,b_c\}\in\mathcal{M}_b$ be a critical bottleneck edge of $\mathcal{G}_b$ relative to $\mathcal{M}_b$.
\end{assum}

%From Assumption~\ref{as:botedge},  $\phi(\mathcal{G}_b,\mathcal{M}_b)$ is the pruned subset of $\mathcal{E}_b$ at the terminating iteration of pruneBAP.

\begin{prop} \label{prop:cluster}
	Consider Assumption~\ref{as:botedge}. The length-one path $P=\{e_c\}$ is a unique alternating path in $\phi(\mathcal{G}_b,\mathcal{M}_b)$ relative to $\mathcal{M}_b$ between $a_c$ and $b_c$.
\end{prop}
\begin{pf}
	Path $P=\{e_c\}$ is trivially an alternating path relative to $\mathcal{M}_b$; edge $e_c\in \mathcal{M}_b$, so $P\subseteq \phi(\mathcal{G}_b,\mathcal{M}_b)$. It remains to show that there does not exist another. Assume for contradiction that there exists an alternating path $P'\neq P$, $P'\subseteq \phi(\mathcal{G}_b,\mathcal{M}_b)$ relative to $\mathcal{M}_b$ between $a_c$ and $b_c$. It follows that $e_c\notin P'$ and therefore $P' \subseteq \phi(\mathcal{G}_b,\mathcal{M}_b)\backslash\{e_c\}$. Furthermore, $P'$ is an augmenting path relative to $\mathcal{M}_b\backslash \{e_c\}$. By Definition~\ref{def:crit}, $e_c$ is not a critical bottleneck edge of $\mathcal{G}_b$, which contradicts Assumption~\ref{as:botedge}. \qed
\end{pf}

The following corollary describes the structure of a bottleneck cluster $\mathcal{G}_b$ based on Proposition~\ref{prop:cluster}.

\begin{cor} \label{cor:tree}
	Consider Assumption~\ref{as:botedge} and let $\mathcal{G}_b$ be a bottleneck cluster with respect to the critical bottleneck edge $e_c$. We form two subgraphs of $\mathcal{G}_b$, denoted as $\mathcal{S}_\mu(\mathcal{G}_b)=(\mathcal{V}_{\mu},\mathcal{E}_{\mu})$ and $\mathcal{S}_\nu(\mathcal{G}_{b})=(\mathcal{V}_{\nu},\mathcal{E}_{\nu})$. Let $\mathcal{V}_{\mu}$ contain the bottleneck agent $a_c$, and let $\mathcal{V}_{\nu}$ contain the bottleneck task $b_c$. Let $\mathcal{V}_b=\mathcal{V}_\mu\cup\mathcal{V}_\nu$ and $\mathcal{V}_\mu\cap\mathcal{V}_\nu=\emptyset$. By Definition~\ref{def:cluster}, it must be possible to construct both $\mathcal{S}_\mu(\mathcal{G}_b)$ and $\mathcal{S}_\nu(\mathcal{G}_b)$ to be alternating trees such that $\mathcal{E}_\mu\cup\mathcal{E}_\nu \subseteq \phi(\mathcal{G}_b,\mathcal{M}_b)$. By Proposition~\ref{prop:cluster}, for all agents $a'\in\mathcal{V}_\mu\cap \mathcal{A}_b$ and for all tasks $b'\in\mathcal{V}_\nu\cap \mathcal{B}_b$, $\{a',b'\}\notin \phi(\mathcal{G}_b,\mathcal{M}_b)$.
\end{cor}

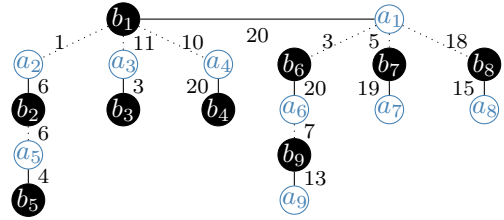
\begin{figure}[thpb]
	\centering
	\begin{tikzpicture}[-,>=stealth',level/.style={sibling distance = 2.5cm/#1,
		level distance = 1.2cm},scale=0.5]
	\node [task] at (0,0) (b1) {$b_1$}
	child{  node [agent] (a2) {$a_2$} edge from parent[dotted]
		child[solid]{ node [task] (b2) {$b_2$} edge from parent[solid]
			child{  node [agent] (a5) {$a_5$} edge from parent[dotted]
				child[solid]{ node [task] (b5) {$b_5$} edge from parent[solid] node[right] {\small$4$}
				} node[right] {\small$6$}
			} node[right] {\small$6$}
		} node[left] {\small$1$}
	}
	child{ node [agent] (a3) {$a_3$} edge from parent[dotted]
		child[solid]{ node [task] (b3) {$b_3$} edge from parent[solid] node[right] {\small$3$}
		} node[right] {\small$11$}
	}
	child{ node [agent] (a4) {$a_4$} edge from parent[dotted]
		child[solid]{ node [task] (b4) {$b_4$} edge from parent[solid] node[left] {\small$20$}
		} node[right] {\small$10$}
	}
	;
	
	\node [agent] at (7,0) (a1) {$a_1$}
	child{  node [task] (b6) {$b_6$} edge from parent[dotted]
		child[solid]{ node [agent] (a6) {$a_6$} edge from parent[solid]
			child{  node [task] (b9) {$b_9$} edge from parent[dotted]
				child[solid]{ node [agent] (a9) {$a_9$} edge from parent[solid] node[right] {\small$13$}
				} node[right] {\small$7$}
			} node[right] {\small$20$}
		} node[left] {\small$3$}
	}
	child{ node [task] (b7) {$b_7$} edge from parent[dotted]
		child[solid]{ node [agent] (a7) {$a_7$} edge from parent[solid] node[left] {\small$19$}
		} node[left] {\small$5$}
	}
	child{ node [task] (b8) {$b_8$} edge from parent[dotted]
		child[solid]{ node [agent] (a8) {$a_8$} edge from parent[solid] node[left] {\small$15$}
		} node[right] {\small$18$}
	}
	;
	\path [draw] (b1) -- node [scale=0.01,label=below:{\small$20$}] {} (a1);
	
	\end{tikzpicture}
	\caption{A bottleneck cluster. Dotted lines represent edges not in matching $\mathcal{M}_b$, solid lines represent edges in $\mathcal{M}_b$. Shown here, a set of agents $\{a_1,a_2,...,a_9\}$ and a set of tasks $\{b_1,b_2,...,b_9\}$. Edge $\{a_1,b_1\}$ is a critical bottleneck edge so Corollary~\ref{cor:tree} applies.}
	\label{fig:trees}
\end{figure}

Fig.~\ref{fig:trees} illustrates Corollary~\ref{cor:tree}. The graph $\mathcal{G}_b$ is represented by two alternating trees $\mathcal{S}_\nu(\mathcal{G}_b)$ and $\mathcal{S}_\mu(\mathcal{G}_b)$ with roots $b_1$ and $a_1$ respectively. The roots $a_1$ and $b_1$ are incident to the critical bottleneck edge $e_c$ with weight $w(e_c)=20$. All edges in both trees are elements of $\phi(\mathcal{G}_b,\mathcal{M}_b)$ as their weights are smaller than or equal to $w(e_c)$ and all edges not in the matching are strictly smaller than $w(e_c)$. Recall from Theorem~\ref{thm:bound} that $w(e_3)\leq \max\{w(e_1),w(e_2)\}$. Given Assumptions~\ref{as:bot1} and~\ref{as:bot2}, the contrapositive of the following lemma provides conditions for the bottleneck edge of graph $\mathcal{G}_3$ to have equal weight to $\max\{w(e_1),w(e_2)\}$.

\begin{lem} \label{lem:reduce}
	Given Assumptions~\ref{as:bot1} and~\ref{as:bot2}, assume both $\mathcal{G}_1$ and $\mathcal{G}_2$ are bottleneck clusters with respect to $e_1$ and $e_2$ respectively. Assume $e_1$ is a critical bottleneck edge of $\mathcal{G}_1$ relative to $\mathcal{M}_1$ and $e_2$ is a critical bottleneck edge of $\mathcal{G}_2$ relative to $\mathcal{M}_2$. Let $w(e_1)\geq w(e_2)$. If $w(e_3)<w(e_1)$, then there exists vertices $i,j\in\mathcal{V}_2$ such that
	\begin{enumerate}[i.]
		\item there exists an edge in $\mathcal{E}_3$ with weight less than $w(e_1)$ between agent $i$ and a task $b'$ in the vertex set of subgraph $\mathcal{S}_\nu(\mathcal{G}_1)$, and
		\item there exists an edge in $\mathcal{E}_3$ with weight less than $w(e_1)$ between task $j$ and an agent $a'$ in the vertex set of subgraph $\mathcal{S}_\mu(\mathcal{G}_1)$, and
		\item there exists an alternating path $P$ between $i$ and $j$ containing only edges with weight less than $w(e_1)$, and $|P\cap\mathcal{M}_2|>|P\backslash\mathcal{M}_2|$.
	\end{enumerate}
\end{lem}

\begin{pf}
	Without loss of generality, let $e_1=\{a_1,b_1\}$. By Proposition~\ref{prop:cluster}, $e_1$ is the only alternating path between $a_1$ and $b_1$ in $\phi(\mathcal{G}_1,\mathcal{M}_1)$. Assume there does not exist vertices $i$ and $j$ such that all i., ii., and iii. are true. Thus, $e_1$ is the only alternating path between $a_1$ and $b_1$ in $\phi(\mathcal{G}_3,\mathcal{M}_1\cup\mathcal{M}_2)$. By Definition~\ref{def:crit}, $e_1$ is also a critical bottleneck edge of $\mathcal{G}_3$ since $e_1\in\arg\max_{e\in\mathcal{M}_1\cup\mathcal{M}_2} w(e)$ and there does not exist an augmenting path in $\phi(\mathcal{G}_3,\mathcal{M}_1\cup\mathcal{M}_2)\backslash\{e_1\}$ relative to $(\mathcal{M}_1\cup\mathcal{M}_2)\backslash\{e_1\}$. Thus, $w(e_3)=w(e_1)$. \qed
\end{pf}

In general, the converse of Lemma~\ref{lem:reduce} does not hold unless we apply some additional assumptions. This leads to the following theorem.

\begin{thm} \label{thm:reduce2}
	Given Assumptions~\ref{as:bot1} and~\ref{as:bot2}, assume both $\mathcal{G}_1$ and $\mathcal{G}_2$ are bottleneck clusters with respect to $e_1$ and $e_2$ respectively. Assume $e_1$ is a critical bottleneck edge of $\mathcal{G}_1$ relative to $\mathcal{M}_1$ and $e_2$ is a critical bottleneck edge of $\mathcal{G}_2$ relative to $\mathcal{M}_2$. Assume that $w(e_1)> w(e_2)$. Let $\arg\max_{\{i,j\}\in\mathcal{M}_1}w(\{i,j\})$ be a singleton. It holds that $w(e_3)<w(e_1)$ if and only if there exists vertices $i,j\in\mathcal{V}_2$ such that conditions i., ii., and iii. from Lemma~\ref{lem:reduce} are true.
\end{thm}
\begin{pf}
	The necessary condition for $w(e_3)<w(e_1)$ holds from Lemma~\ref{lem:reduce}. We now prove the sufficient condition. Assume there exists vertices $i$ and $j$ such that all i., ii., and iii. are true. Then, aside from path $P=\{e_1=\{a_1,b_1\}\}$, there exists another alternating path $P'$ between $a_1$ and $b_1$, which does not contain the edge $e_1$. Namely, the alternating path $P'$ constructed from the union of the alternating paths between $a_1$ and $b'$, $b'$ and $i$, $i$ and $j$, $j$ and $a'$, and $a'$ and $b_1$. Thus, there exists an augmenting path $P'\subseteq\phi(\mathcal{G}_3,\mathcal{M}_1\cup\mathcal{M}_2)\backslash\{e_1\}$ relative to $(\mathcal{M}_1\cup\mathcal{M}_2)\backslash\{e_1\}$. From Lemma~\ref{lem:distBAP}, there exists an MCM $\mathcal{M}'$ of $\mathcal{G}_3$ such that $\mathcal{M}'\in \phi(\mathcal{G}_3,\mathcal{M}_1\cup\mathcal{M}_2)\backslash\{e_1\}$. By the assumptions on $e_1$, $\phi(\mathcal{G}_3,\mathcal{M}_1\cup\mathcal{M}_2)\backslash\{e_1\}$ contains only edges with weights strictly smaller than $w(e_1)$. Thus, there exists an MCM of $\mathcal{G}_3$ with all edges having weight smaller than $w(e_1)$, i.e., $w(e_3)$ must be smaller than $w(e_3)$. \qed
\end{pf}

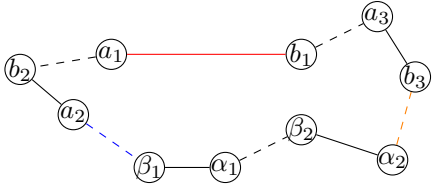
\begin{figure}[thpb]
	\centering
	\begin{tikzpicture}
	\node[draw,circle,minimum size=0.4cm] at (-0.5,0) (a1) {};
	\node[draw,circle,minimum size=0.4cm] at (2,0) (b1) {};
	
	\node[draw,circle,minimum size=0.4cm] at (-1,-0.8) (a2) {};
	\node[draw,circle,minimum size=0.4cm] at (-1.7,-0.2) (b2) {};
	
	\node[draw,circle,minimum size=0.4cm] at (3,0.5) (a3) {};
	\node[draw,circle,minimum size=0.4cm] at (3.5,-0.3) (b3) {};
	
	\node[draw,circle,minimum size=0.4cm] at (1,-1.5) (aa1) {};
	\node[draw,circle,minimum size=0.4cm] at (0,-1.5) (bb1) {};
	
	\node[draw,circle,minimum size=0.4cm] at (3.2,-1.4) (aa2) {};
	\node[draw,circle,minimum size=0.4cm] at (2,-1) (bb2) {};
	
	\node[] at (a1) {$a_1$};
	\node[] at (2,0) {$b_1$};
	
	\node[] at (a2) {$a_2$};
	\node[] at (b2) {$b_2$};
	
	\node[] at (a3) {$a_3$};
	\node[] at (b3) {$b_3$};
	
	\node[] at (aa1) {$\alpha_1$};
	\node[] at (bb1) {$\beta_1$};
	
	\node[] at (aa2) {$\alpha_2$};
	\node[] at (bb2) {$\beta_2$};
	
	\path[draw,red] (a1) -- (b1);
	\path[draw] (a2) -- (b2);
	\path[draw] (a3) -- (b3);
	\path[draw] (aa1) -- (bb1);
	\path[draw] (aa2) -- (bb2);
	
	\path[draw,dashed] (a1) -- (b2);
	\path[draw,dashed] (b1) -- (a3);
	\path[draw,dashed] (aa1) -- (bb2);
	
	\path[draw,dashed,blue] (a2) -- (bb1);
	\path[draw,dashed,orange] (aa2) -- (b3);
	\end{tikzpicture}
	\caption{An illustration of Theorem~\ref{thm:reduce2}. We have $\mathcal{V}_1=\{a_1,a_2,a_3\}\cup\{b_1,b_2,b_3\}$ and $\mathcal{G}_1$ is a bottleneck cluster with respect to $e_1=\{a_1,b_1\}$. Edge $e_1$ is shown as a solid red line. We have $\mathcal{V}_2=\{\alpha_1,\alpha_2\}\cup\{\beta_1,\beta_2\}$ and $\mathcal{G}_2$ is a bottleneck cluster with respect to $e_2=\{\alpha_2,\beta_2\}$. The length of each line corresponds to the weight of that edge. Solid lines show edges in $\mathcal{M}_1\cup\mathcal{M}_2$. Dashed lines show edges in $\mathcal{M}_3$ and $w(e_3)<w(e_1)$.}
	\label{fig:warmstart}
\end{figure}

Fig.~\ref{fig:warmstart} illustrates the sufficient condition of Theorem~\ref{thm:reduce2}. The orange dashed line is an edge that satisfies the condition i.\ since $\alpha_2\in \mathcal{V}_2$ and there exists edge $\{\alpha_2,b_3\}\in \mathcal{V}_3$, where $b_3$ is a task in the vertex set of $\mathcal{S}_\nu(\mathcal{G}_1)$. The blue dashed line satisfies condition ii.\ since $\beta_1\in \mathcal{V}_2$ and there exists edge $\{\beta_1,a_2\}\in \mathcal{V}_3$, where $a_2$ is an agent in the vertex set of $\mathcal{S}_\mu(\mathcal{G}_1)$. Condition iii. is satisfied since there is an alternating path $P=\{\{\alpha_1,\beta_1\},\{\alpha_2,\beta_2\}, \{\alpha_1,\beta_2\}\}$ between $\beta_1$ and $\alpha_2$, and $|\{\{\alpha_1,\beta_1\},\{\alpha_2,\beta_2\}\}|>|\{\{\alpha_1,\beta_2\}\}|$, i.e., $P$ starts with a dashed line and ends with a dashed line. Corollary~\ref{cor:main} follows from Theorem~\ref{thm:bound} and~\ref{thm:reduce2}.

\begin{cor} \label{cor:main}
	If one or more of conditions i., ii., or iii. in Theorem~\ref{thm:reduce2} do not hold, then $\mathcal{M}_3=\mathcal{M}_1\cup\mathcal{M}_2$ is a solution to Problem~\ref{prob:sol1}. 
\end{cor}

\section{Algorithm for Solving the BAP} \label{section:distBAP}

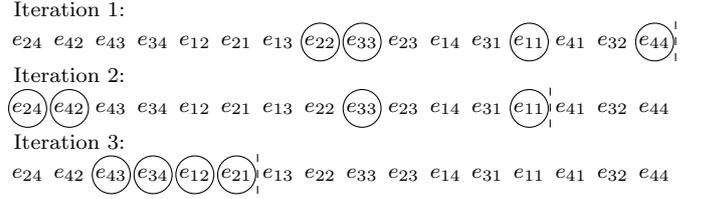
\begin{figure}[thpb]
	\centering
	\begin{tikzpicture}[scale=0.55]
		\tikzstyle{every node}=[font=\small]
		\node [] at (1,0.8) {Iteration 1:};
		\node [] at (0,0) (w1) {$e_{24}$};
		\node [] at (1,0) (w2) {$e_{42}$};
		\node [] at (2,0) (w3) {$e_{43}$};
		\node [] at (3,0) (w4) {$e_{34}$};
		\node [] at (4,0) (w5) {$e_{12}$};
		\node [] at (5,0) (w6) {$e_{21}$};
		\node [] at (6,0) (w7) {$e_{13}$};
		\node [] at (7,0) (w8) {$e_{22}$};
		\node [] at (8,0) (w9) {$e_{33}$};
		\node [] at (9,0) (w10) {$e_{23}$};
		\node [] at (10,0) (w11) {$e_{14}$};
		\node [] at (11,0) (w12) {$e_{31}$};
		\node [] at (12,0) (w13) {$e_{11}$};
		\node [] at (13,0) (w14) {$e_{41}$};
		\node [] at (14,0) (w15) {$e_{32}$};
		\node [] at (15,0) {$e_{44}$};
		
		\node [draw,circle,minimum size = 0.5cm] at (15,0) {};
		\node [draw,circle,minimum size = 0.5cm] at (12,0) {};
		\node [draw,circle,minimum size = 0.5cm] at (8,0) {};
		\node [draw,circle,minimum size = 0.5cm] at (7,0) {};
		\draw [dashed] (15.5,0.5) -- (15.5,-0.5);
		%-------------
		\node [] at (1,-0.8) {Iteration 2:};
		\node [] at (0,-1.6) (w1) {$e_{24}$};
		\node [] at (1,-1.6) (w2) {$e_{42}$};
		\node [] at (2,-1.6) (w3) {$e_{43}$};
		\node [] at (3,-1.6) (w4) {$e_{34}$};
		\node [] at (4,-1.6) (w5) {$e_{12}$};
		\node [] at (5,-1.6) (w6) {$e_{21}$};
		\node [] at (6,-1.6) (w7) {$e_{13}$};
		\node [] at (7,-1.6) (w8) {$e_{22}$};
		\node [] at (8,-1.6) (w9) {$e_{33}$};
		\node [] at (9,-1.6) (w10) {$e_{23}$};
		\node [] at (10,-1.6) (w11) {$e_{14}$};
		\node [] at (11,-1.6) (w12) {$e_{31}$};
		\node [] at (12,-1.6) (w13) {$e_{11}$};
		\node [] at (13,-1.6) (w14) {$e_{41}$};
		\node [] at (14,-1.6) (w15) {$e_{32}$};
		\node [] at (15,-1.6) {$e_{44}$};
		
		\node [draw,circle,minimum size = 0.5cm] at (12,-1.6) {};
		\node [draw,circle,minimum size = 0.5cm] at (0,-1.6) {};
		\node [draw,circle,minimum size = 0.5cm] at (8,-1.6) {};
		\node [draw,circle,minimum size = 0.5cm] at (1,-1.6) {};
		\draw [dashed] (12.5,-1.1) -- (12.5,-2.1);
		
		%------------
		\node [] at (1,-2.4) {Iteration 3:};
		\node [] at (0,-3.2) (w1) {$e_{24}$};
		\node [] at (1,-3.2) (w2) {$e_{42}$};
		\node [] at (2,-3.2) (w3) {$e_{43}$};
		\node [] at (3,-3.2) (w4) {$e_{34}$};
		\node [] at (4,-3.2) (w5) {$e_{12}$};
		\node [] at (5,-3.2) (w6) {$e_{21}$};
		\node [] at (6,-3.2) (w7) {$e_{13}$};
		\node [] at (7,-3.2) (w8) {$e_{22}$};
		\node [] at (8,-3.2) (w9) {$e_{33}$};
		\node [] at (9,-3.2) (w10) {$e_{23}$};
		\node [] at (10,-3.2) (w11) {$e_{14}$};
		\node [] at (11,-3.2) (w12) {$e_{31}$};
		\node [] at (12,-3.2) (w13) {$e_{11}$};
		\node [] at (13,-3.2) (w14) {$e_{41}$};
		\node [] at (14,-3.2) (w15) {$e_{32}$};
		\node [] at (15,-3.2) {$e_{44}$};
		
		\node [draw,circle,minimum size = 0.5cm] at (4,-3.2) {};
		\node [draw,circle,minimum size = 0.5cm] at (5,-3.2) {};
		\node [draw,circle,minimum size = 0.5cm] at (3,-3.2) {};
		\node [draw,circle,minimum size = 0.5cm] at (2,-3.2) {};
		\draw [dashed] (5.5,-2.7) -- (5.5,-3.7);
	\end{tikzpicture}
	\caption{A demonstration of pruneBAP with $\mathcal{A}_3=\{a_1,a_2,a_3,a_4\}$ and $\mathcal{B}_3=\{b_1,b_2,b_3,b_4\}$. Edges in $\mathcal{E}_3$ are arranged in order of ascending weight, where $e_{pq}$ is the edge between agent $a_p$ and task $b_q$. At iteration 1, the initial arbitrary MCM is denoted by the 4 circled edges. Edges to the right of the dashed lines have been pruned from $\mathcal{E}_3$. Note, $w(e_{44})\geq w(e_{11})\geq w(e_{21})$, i.e., with each iteration the weight of the largest edge in the current MCM is non-increasing. The algorithm terminates when a matching of size 4 does not exist in the remaining edges to the left of the dashed line.}
	\label{fig:distBAP}
\end{figure}

In this section, we discuss how the algorithm from \cite{distBAP} makes use of Assumptions~\ref{as:bot1} and~\ref{as:bot2} to solve Problem~\ref{prob:sol1}. Let us refer to this algorithm as pruneBAP. Fig.~\ref{fig:distBAP} is an illustration of this algorithm.
%The pruneBAP algorithm solves $BOT(\mathcal{G}_3)$ and we first summarise pruneBAP in three steps. The algorithm is initialised with an arbitrary initial MCM $\mathcal{M}_0$. Then, the following steps are repeated until convergence to a solution of $BOT(\mathcal{G}_3)$.
%\begin{enumerate}[Step 1]
%	\item Prune the edge set $\mathcal{E}_3$ of $\mathcal{G}_3$ by removing edges; given an MCM $\mathcal{M}$, obtain the set $\phi(\mathcal{G}_3,\mathcal{M})\in \mathcal{E}_3$.
%	\item Remove one edge $e_l$ with largest weight from $\phi(\mathcal{G}_3,\mathcal{M})$. Based on Step 1, this edge must also be an element of the current MCM $\mathcal{M}$.
%	\item Find a new MCM of $\mathcal{G}_3$ from the remaining edges in $\phi(\mathcal{G}_3,\mathcal{M})\backslash\{e_l\}$, by making use of an augmenting path and $\mathcal{M}\backslash\{e_l\}$. See Lemma~\ref{lem:distBAP} for details.
%\end{enumerate}
%
%When Step 3 does not yield an MCM, the algorithm terminates. The last valid MCM is a solution to $BOT(\mathcal{G}_3)$.

%\begin{rem}
%	Note that the steps in pruneBAP are amenable for constructing a bottleneck cluster. In fact, the subgraph $\mathcal{S}_\nu(\mathcal{G}_3)$ discussed in Corollary~\ref{cor:tree} is constructed as an artefact of implementing pruneBAP. Then, to obtain a bottleneck cluster it remains to construct $\mathcal{S}_\mu$, which involves one additional graph search to find all vertices with an alternating path to the bottleneck agent $a_l$.
%\end{rem}

\subsection{Warm-starting Versus Cold-starting pruneBAP}

\begin{figure}[thpb]
	\centering
	\begin{tikzpicture}[scale=0.55]
	\tikzstyle{every node}=[font=\small]
	\node [] at (2.5,0.8) {$\mathcal{M}_1$ from Assumption \ref{as:bot1}:};
	\node [] at (4,0) (w2) {$e_{12}$};
	\node [] at (5,0) (w3) {$e_{21}$};
	\node [] at (7,0) (w8) {$e_{22}$};
	\node [] at (12,0) (w9) {$e_{11}$};
	
	\node [draw,circle,minimum size = 0.5cm] at (4,0) {};
	\node [draw,circle,minimum size = 0.5cm] at (5,0) {};
	\draw [dashed] (5.5,0.5) -- (5.5,-0.5);
	%-------------
	\node [] at (2.5,-0.8) {$\mathcal{M}_2$ from Assumption \ref{as:bot2}:};
	\node [] at (2,-1.6) (w1) {$e_{43}$};
	\node [] at (3,-1.6) (w4) {$e_{34}$};
	\node [] at (8,-1.6) (w6) {$e_{33}$};
	\node [] at (15,-1.6) (w7) {$e_{44}$};
	
	\node [draw,circle,minimum size = 0.5cm] at (2,-1.6) {};
	\node [draw,circle,minimum size = 0.5cm] at (3,-1.6) {};
	\draw [dashed] (3.5,-1.1) -- (3.5,-2.1);

	%------------
	\node [] at (2.5,-2.4) {Warm-start, Iteration 1:};
	\node [] at (0,-3.2) (w1) {$e_{24}$};
	\node [] at (1,-3.2) (w2) {$e_{42}$};
	\node [] at (2,-3.2) (w3) {$e_{43}$};
	\node [] at (3,-3.2) (w4) {$e_{34}$};
	\node [] at (4,-3.2) (w5) {$e_{12}$};
	\node [] at (5,-3.2) (w6) {$e_{21}$};
	\node [] at (6,-3.2) (w7) {$e_{13}$};
	\node [] at (7,-3.2) (w8) {$e_{22}$};
	\node [] at (8,-3.2) (w9) {$e_{33}$};
	\node [] at (9,-3.2) (w10) {$e_{23}$};
	\node [] at (10,-3.2) (w11) {$e_{14}$};
	\node [] at (11,-3.2) (w12) {$e_{31}$};
	\node [] at (12,-3.2) (w13) {$e_{11}$};
	\node [] at (13,-3.2) (w14) {$e_{41}$};
	\node [] at (14,-3.2) (w15) {$e_{32}$};
	\node [] at (15,-3.2) {$e_{44}$};
	
	\node [draw,circle,minimum size = 0.5cm] at (2,-3.2) {};
	\node [draw,circle,minimum size = 0.5cm] at (3,-3.2) {};
	\node [draw,circle,minimum size = 0.5cm] at (4,-3.2) {};
	\node [draw,circle,minimum size = 0.5cm] at (5,-3.2) {};
	\draw [dashed] (5.5,-2.7) -- (5.5,-3.7);
	\end{tikzpicture}
	\caption{A demonstration a warm-start. Here, $\mathcal{A}_1=\{a_1,a_2\}$, $\mathcal{A}_2=\{a_3,a_4\}$, $\mathcal{B}_1=\{b_1,b_2\}$, and $\mathcal{B}_2=\{b_3,b_4\}$. Then, $\mathcal{A}_3$ and $\mathcal{B}_3$ are the same as in Fig.~\ref{fig:distBAP}. In this example, warm-starting pruneBAP allows the solution to $BOT(\mathcal{G}_3)$ to be found in 1 iteration.}
	\label{fig:combine}
\end{figure}
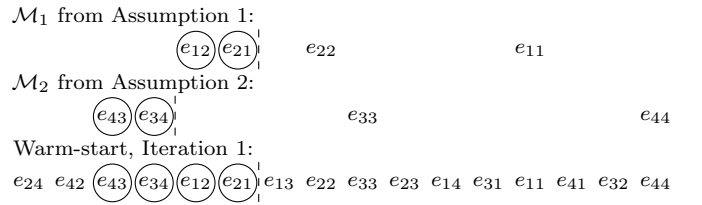

Solving Problem~\ref{prob:sol1} by pruneBAP with an arbitrary MCM $\mathcal{M}_0$ at initialisation does not make use of Assumptions \ref{as:bot1} and \ref{as:bot2}. We denote this as a cold-start to pruneBAP. Given Assumptions \ref{as:bot1} and \ref{as:bot2}, consider the following. It holds that the set $\tilde{\mathcal{M}}:=\mathcal{M}_1\cup\mathcal{M}_2$ is an MCM of the graph $\mathcal{G}_3$. Without loss of generality, let $w(e_1)\geq w(e_2)$. Then, it also holds that $e_1$ is the largest edge in $\tilde{\mathcal{M}}$. We use make use of $\tilde{\mathcal{M}}$ to solve Problem 1 by choosing it as the initial MCM of pruneBAP. Edges in the set $\{e\in\mathcal{E}_3|w(e)\geq w(e_1),e\notin \mathcal{M}_0\}$ are removed from $\mathcal{E}_3$ in the first iteration. We denote this as a warm-start to pruneBAP. Fig.~\ref{fig:combine} illustrates a warm-start to pruneBAP.

\begin{rem}
	%Let $\mathcal{M}_0=\mathcal{M}\in\mathcal{C}(\mathcal{G}_3)$ denote an arbitrary MCM used to cold-start pruneBAP. Under Assumptions \ref{as:bot1} and \ref{as:bot2}, let $\tilde{\mathcal{M}}=\mathcal{M}_1\cup\mathcal{M}_2$ denote the MCM used to warm-start pruneBAP. Assume $w(e_1)\geq w(e_2)$. There is no guarantee that $\max\limits_{(i,j)\in\mathcal{M}_0}w(\{i,j\})>w(e_1)$. Thus, there is no guarantee that a warm-start allows the removal of more edges from $\mathcal{E}_3$ at iteration 1.
	Warm-starting is a heuristic, a warm-start does not guarantee fewer iterations for convergence to a solution to $BOT(G\mathcal{G}_3)$. For a cold-start, we choose an arbitrary initial MCM $\mathcal{M}_{cold}$; by chance this $\mathcal{M}_{cold}$ could be the solution to $BOT(\mathcal{G}_3)$.
\end{rem}

Given Assumptions~\ref{as:bot1} and~\ref{as:bot2}, warm-starting is a way to make use of the available information to solve Problem~\ref{prob:sol1}.

\section{Case Studies}

Consider agents and tasks represented by points in a vector space $S$ with a distance function $D:S\times S\mapsto \mathbb{R}^+$. For example, this could be a ride-sharing application, where agents are vehicles and their tasks are to pick up customers. Here we consider a 2-dimensional space  $S=\mathbb{R}^2$ and Euclidean distance $D(x,y)=\lVert x-y \rVert_2$. Agents are to be assigned to move from their initial positions to target destinations based on the BAP with distance as weights.

\subsection{Task Reassignment}

\begin{casest} \label{sc:random}
	Let $\mathcal{A}=\{a_1,a_2,...,a_{m_3}\}\subset S$ be the initial locations of a set of agents. Let $\mathcal{B}_1=\{b_1,b_2,...,b_{n_1}\}\subset S$ be the set of goal locations. Assume $m_3>n_1$. We first solve $BOT((\mathcal{A}\cup\mathcal{B}_1,\mathcal{E}))$, where $\mathcal{E}=\{\{i,j\}|i\in\mathcal{A},j\in\mathcal{B}_1\}$ to determine an assignment of tasks to agents that minimises the worst-case distance an agent must travel to reach a goal location. Without loss of generality, assume vehicles at positions $\mathcal{A}_1=\{a_1,a_2,...,a_{n_1}\}$ are assigned to goals at $\mathcal{B}_1$. Now assume that a second set of goal locations becomes available to agents. Let $\mathcal{B}_2=\{\beta_1,\beta_2,...,\beta_{n_2}\}\subset S$ be the set of new goal locations. Assume that $m_3\geq n_1+n_2$. Let $\mathcal{A}_2=\{a_{n_1+1},a_{n_1+2},...,a_{m_3}\}$ be the locations of the remaining unassigned agents. We now assign the new goals to the remaining agents, i.e., solve $BOT((\mathcal{A}_2\cup\mathcal{B}_2,\mathcal{E}_2))$, where $\mathcal{E}_2=\{\{i,j\}|i\in\mathcal{A}_2,j\in\mathcal{B}_2\}$.
\end{casest}

\begin{figure}[thpb]
	\centering
	\includegraphics[scale=0.45]{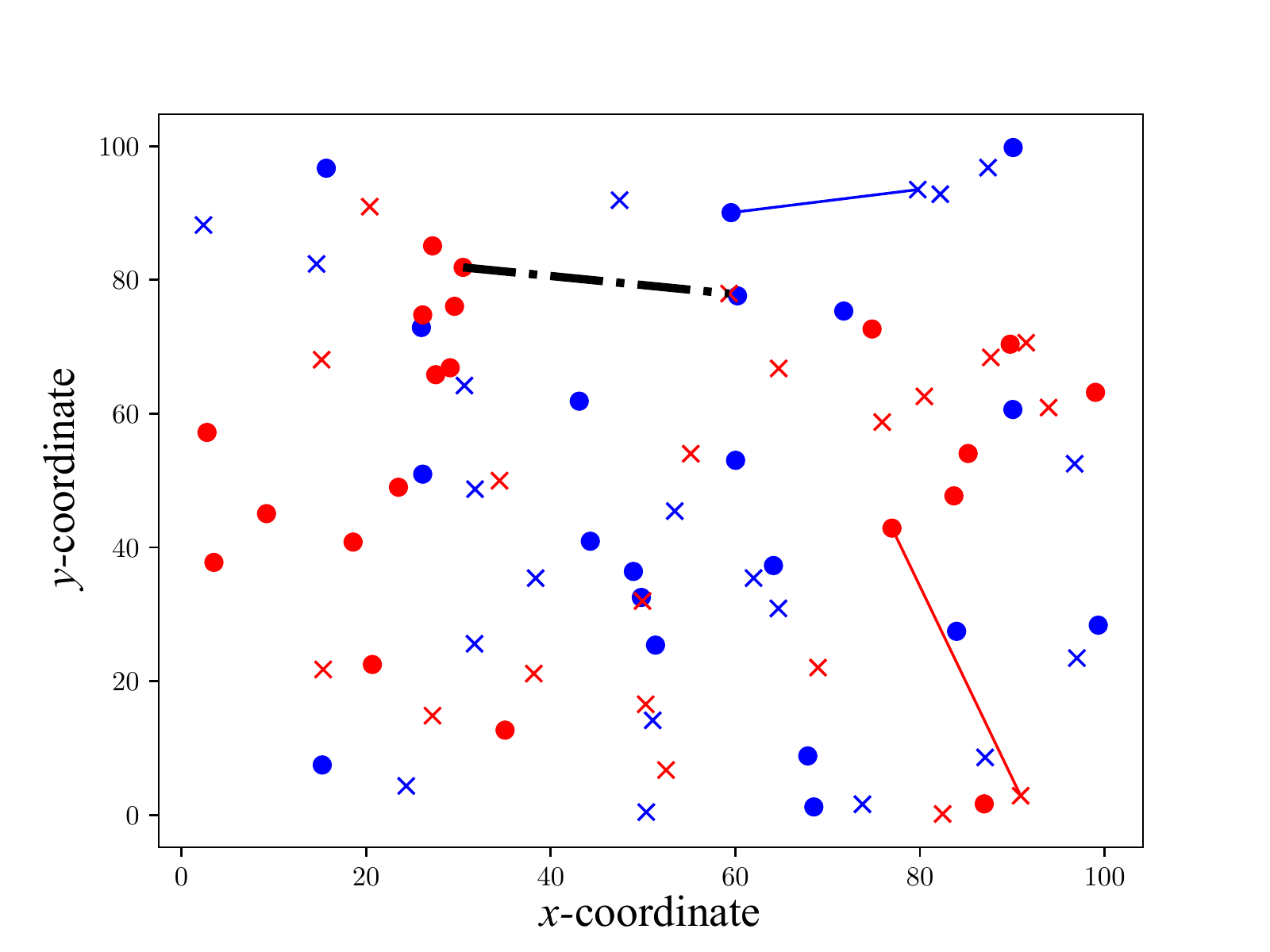}
	\caption{Case Study 1: Sample configuration of agents and tasks. The locations in $\mathcal{A}_1$, $\mathcal{A}_2$, $\mathcal{B}_1$, and $\mathcal{B}_2$ are represented by blue dots, red dots, blue crosses and red crosses respectively. The bottleneck edges from solving $BOT((\mathcal{A}\cup\mathcal{B}_1,\mathcal{E}))$, $BOT((\mathcal{A}_2\cup\mathcal{B}_2,\mathcal{E}_2))$ and $BOT((\mathcal{A}\cup\mathcal{B}_1\cup\mathcal{B}_2,\mathcal{E}_3))$ are shown by the solid blue, solid red and black dash-dotted lines respectively.}
	\label{fig:uniform}
\end{figure}

\begin{figure}[thpb]
	\centering
	\includegraphics[scale=0.45]{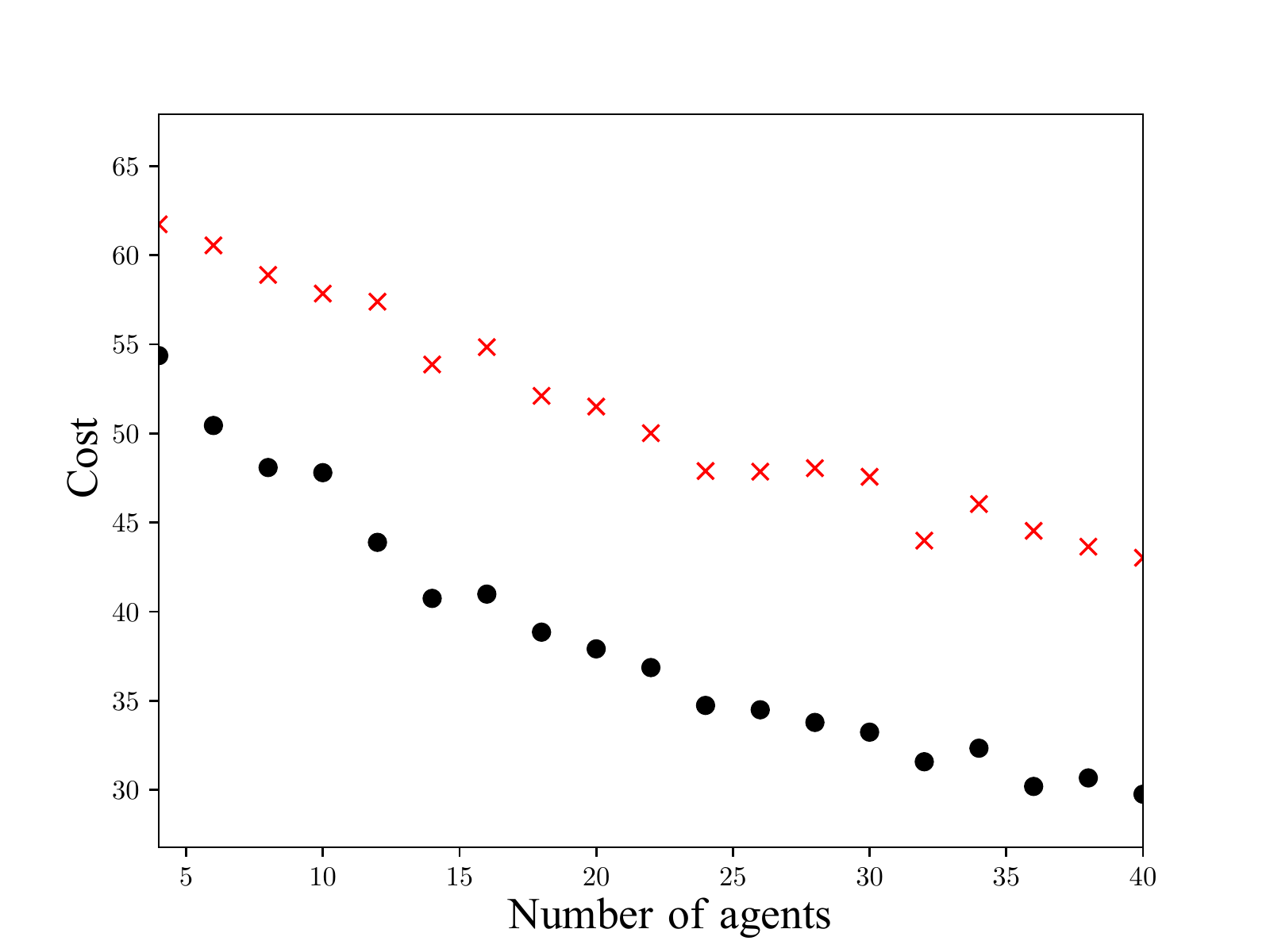}
	\caption{Case Study 1: Cost of assignments. Red crosses show the average weight of the largest edge in the MCM used to warm-start pruneBAP, from 100 simulations. Black dots show the average weight of the bottleneck edge from solving $BOT((\mathcal{A}\cup\mathcal{B}_1\cup\mathcal{B}_2,\mathcal{E}_3))$, from 100 simulations.}
	\label{fig:graphwarm}
\end{figure}

By Theorem~\ref{thm:bound}, the assignment obtained from solving $BOT((\mathcal{A}\cup\mathcal{B}_1,\mathcal{E}))$ and $BOT((\mathcal{A}_2\cup\mathcal{B}_2,\mathcal{E}_2))$ in Case Study~\ref{sc:random} is not necessarily the optimal solution to $BOT((\mathcal{A}\cup\mathcal{B}_1\cup\mathcal{B}_2,\mathcal{E}_3))$, where $\mathcal{E}_3=\{\{i,j\}|i\in\mathcal{A},j\in \mathcal{B}_1\cup\mathcal{B}_2\}$. Fig.~\ref{fig:uniform} shows a numerical example of a case where the optimal assignment is of lower cost than the assignment used to warm-start pruneBAP. For this example, $m_3=40$, $n_1=20$ and $n_2=20$. The data was generated using continuous uniform distributions with range $[0,100)$ for both coordinates $x$ and $y$. Fig.~\ref{fig:graphwarm} shows a plot of the average cost of the assignment used as warm-start to initialise pruneBAP and the average cost of the optimal assignment after pruneBAP has terminated. For all simulations, $m_3=n_1+n_2$. For each even value of $m_3$, 100 simulations were generated. We observe that the cost of the assignment obtained from the subproblems is never greater than the cost of an optimal solution to $BOT(\mathcal{G}_3)$, in accordance with Theorem~\ref{thm:bound}. In this case, the unstructured distribution of the locations results in all of the conditions in Theorem~\ref{thm:reduce2} being satisfied and we observe that $w(e_3)<\max\{w(e_1),w(e_2)\}$ as expected.

\subsection{Clustering of Agents and Tasks}

\begin{casest} \label{sc:cluster}
	Let $\mathcal{A}_1=\{a_1,a_2,...,a_{m_1}\}\subset S$ and $\mathcal{A}_2=\{\alpha_1,\alpha_2,...,\alpha_{m_2}\}\subset S$ be the initial locations of two sets of agents. Let $\mathcal{B}_1=\{b_1,b_2,...,b_{n_1}\}\subset S$ and $\mathcal{B}_2=\{\beta_1,\beta_2,...,\beta_{n_2}\}\subset S$ be the sets of goal locations. Assume that $m_1\geq n_1$ and $m_2\geq n_2$, i.e., there are more agents than there are goals. Assume the set of locations $\mathcal{A}_1$ and $\mathcal{B}_1$ are separated geographically from $\mathcal{A}_2$ and $\mathcal{B}_2$.
\end{casest}

\begin{figure}[thpb]
	\centering
	\includegraphics[scale=0.45]{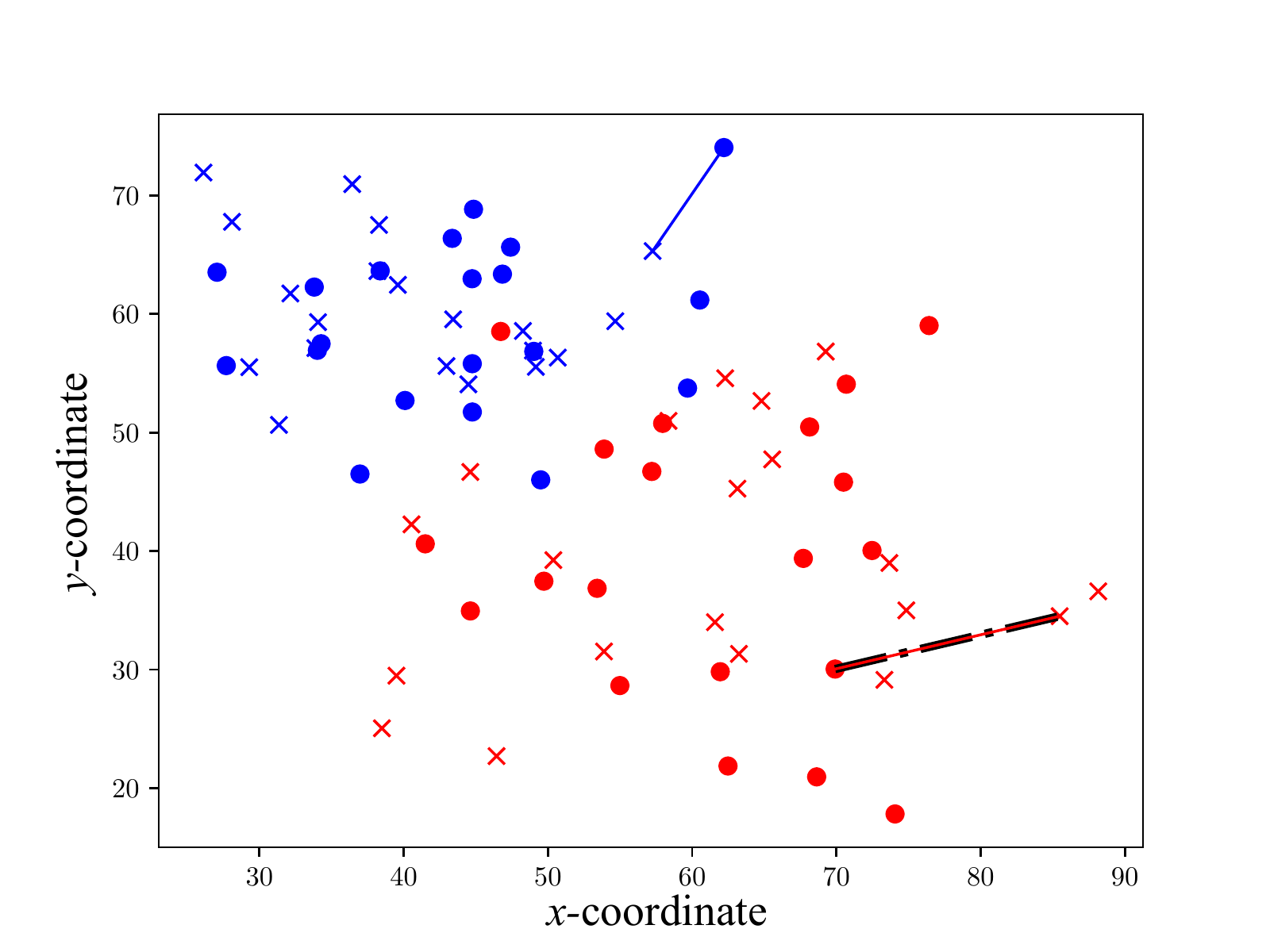}
	\caption{Case Study 2: Sample configuration of agents and tasks. The positions of agents and tasks $\mathcal{A}_1$, $\mathcal{A}_2$, $\mathcal{B}_1$, and $\mathcal{B}_2$ are represented by blue dots, red dots, blue crosses and red crosses respectively. The solid blue line shows the bottleneck edge from solving $BOT(\mathcal{G}_1)$. The solid red line shows the bottleneck edge from solving $BOT(\mathcal{G}_2)$. The black dash-dotted line shows the bottleneck edge from solving $BOT(\mathcal{G}_3)$.}
	\label{fig:cluster}
\end{figure}

\begin{figure}[thpb]
	\centering
	\includegraphics[scale=0.45]{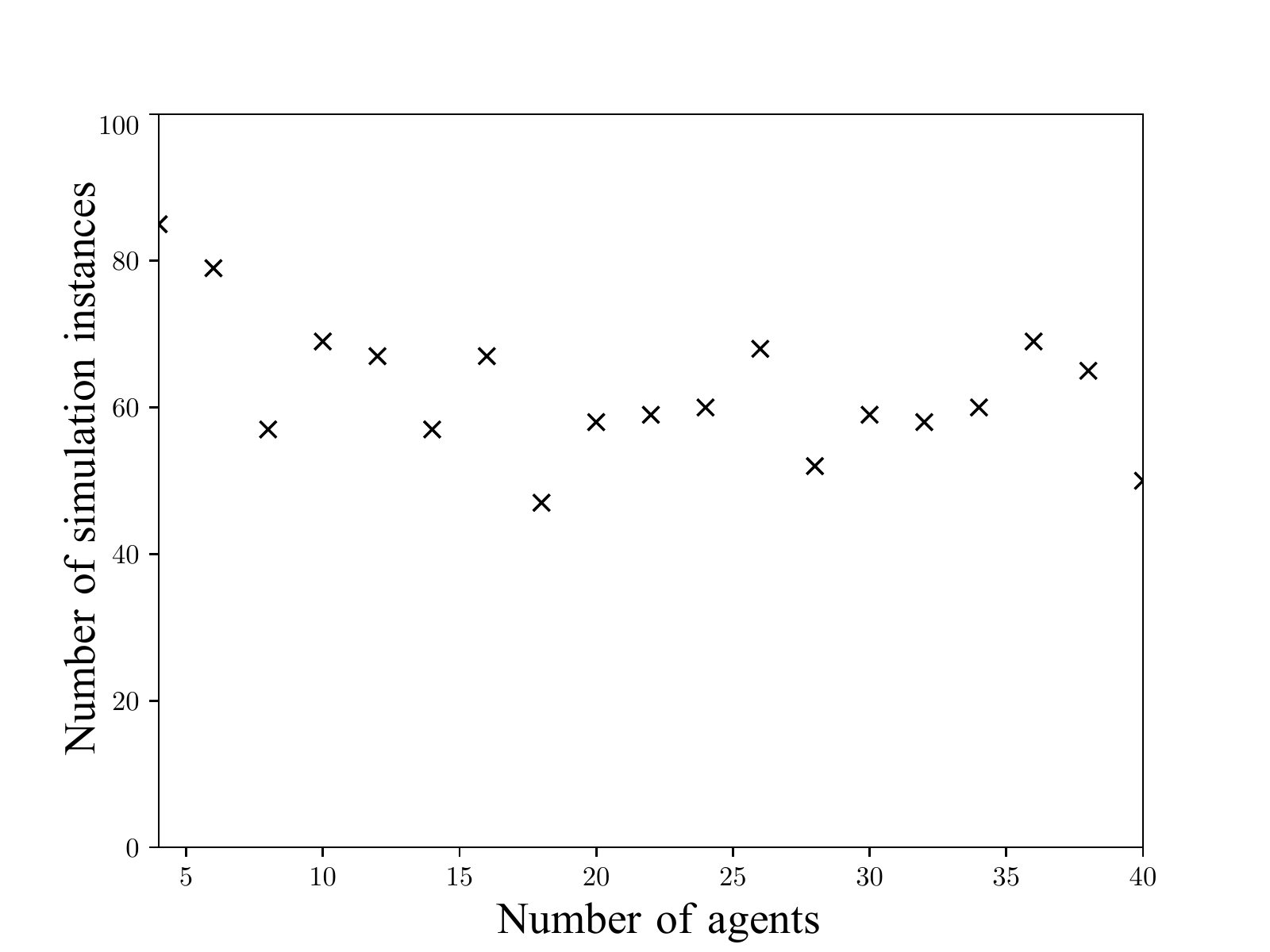}
	\caption{Case Study 2: Empirical probability of solution to subproblems resulting in exact solution of the combined problem. Each cross shows the number of instances out of 100 simulations that the MCM $\mathcal{M}_3=\mathcal{M}_1\cup\mathcal{M}_2$ is an optimal solution to $BOT(\mathcal{G}_3)$. The number of agents represents $m_1+m_2$.}
	\label{fig:reassignment}
\end{figure}

In Case Study~\ref{sc:cluster}, we illustrate an example where not all of the conditions i., ii., and iii. in Lemma~\ref{lem:reduce} hold. Fig.~\ref{fig:cluster} shows a numerical example where the initial assignment $\mathcal{M}_1\cup \mathcal{M}_2$ used to warm-start pruneBAP is in fact the optimal assignment of $\mathcal{B}_1\cup \mathcal{B}_2$ to $\mathcal{A}_1\cup \mathcal{A}_2$. In this example, $m_1=20$, $m_2=20$, $n_1=20$ and $n_2=20$. The data in Fig.~\ref{fig:cluster} was generated using independant normal distributions with a variance of 100 for each distribution. The distributions for sets $\mathcal{A}_1$ and $\mathcal{B}_1$ are centred at the point $(x,y)=(40,60)$. The distributions for sets $\mathcal{A}_2$ and $\mathcal{B}_2$ are centred at the point $(x,y)=(60,40)$. Fig.~\ref{fig:reassignment} shows the number of instances out of 100 simulations for which the behaviour in Fig.~\ref{fig:cluster} is observed. That is, the instances where the bottleneck edges obtained from the subproblems directly results in an optimal solution to $BOT(\mathcal{G}_3)$, where $\mathcal{G}_3$ is defined as in Section~\ref{formulation}. The number of agents equals the number of tasks for each simulation, i.e., $m_1=n_1=m_2=n_2$. For each simulation, positions were generated using the same normal distribution as in Fig.~\ref{fig:cluster}. We now observe realisations where the cost of the assignment obtained from the subproblems is equal to the cost of an optimal solution to $BOT(\mathcal{G}_3)$. This illustrates that with this distribution of agents and tasks there are instances where there is structure such that the conditions in Theorem~\ref{thm:reduce2} do not all hold and $w(e_3)=\max\{w(e_1),w(e_2)\}$.

\section{Conclusion}

We discussed properties of pruneBAP that allow us to warm-start the algorithm given BAP solutions to divided sets of tasks and agents. The solutions based on the divided problems forms an MCM of the combined problem. The pruneBAP algorithm can be initialised with any MCM, and thus allows us to make use of the solutions based on the divided sets. We then have an upper bound on the BAP solution to the combined problem in terms of the bottleneck edges of the divided problems. We also introduced the novel concept of a bottleneck cluster relative to a bottleneck edge. This idea is inspired by the pruneBAP algorithm and the alternating tree that is obtained as a result of the algorithm. Using bottleneck clusters, we provided conditions such that the initial MCM used to warm-start pruneBAP is a solution to the BAP. From numerical simulations motivated by ride-sharing, we illustrate an example where the conditions hold if there exist clusters that are separated in space.

An interesting future direction is the investigation of methods to optimally partition agents and tasks. Another direction would be to investigate clustering properties for the LAP.

\bibliography{ifacconf}             % bib file to produce the bibliography
                                                     % with bibtex (preferred)
                                                   
%\begin{thebibliography}{xx}  % you can also add the bibliography by hand

%\bibitem[Able(1956)]{Abl:56}
%B.C. Able.
%\newblock Nucleic acid content of microscope.
%\newblock \emph{Nature}, 135:\penalty0 7--9, 1956.

%\bibitem[Able et~al.(1954)Able, Tagg, and Rush]{AbTaRu:54}
%B.C. Able, R.A. Tagg, and M.~Rush.
%\newblock Enzyme-catalyzed cellular transanimations.
%\newblock In A.F. Round, editor, \emph{Advances in Enzymology}, volume~2, pages
%  125--247. Academic Press, New York, 3rd edition, 1954.

%\bibitem[Keohane(1958)]{Keo:58}
%\newblock \emph{Power and Interdependence: World Politics in Transitions}.
%\newblock Little, Brown \& Co., Boston, 1958.

%\bibitem[Powers(1985)]{Pow:85}
%\newblock Is there a way out?
%\newblock \emph{Harpers}, pages 35--47, June 1985.

%\bibitem[Soukhanov(1992)]{Heritage:92}
%\newblock \emph{{The American Heritage. Dictionary of the American Language}}.
%\newblock Houghton Mifflin Company, 1992.

%\end{thebibliography}

%\appendix
%\section{A summary of Latin grammar}    % Each appendix must have a short title.
%\section{Some Latin vocabulary}              % Sections and subsections are supported  
                                                                         % in the appendices.
\end{document}